\documentclass[11pt,letterpaper]{amsart}
\usepackage{amscd}
\usepackage{amssymb}
\usepackage{amsthm}
\usepackage{enumerate}
\usepackage[normalem]{ulem}
\usepackage{mathtools}
\usepackage[usenames,dvipsnames]{xcolor}
\usepackage{stmaryrd}
\usepackage[left=1in,top=1in,right=1in]{geometry}
\usepackage[T2A]{fontenc}
\usepackage[utf8]{inputenc}
\usepackage{mathrsfs}
\usepackage[bookmarks=false]{hyperref}

\newtheorem{thm}{Theorem}[section]
\newtheorem{lemm}[thm]{Lemma}
\newtheorem{prop}[thm]{Proposition}
\newtheorem{cor}[thm]{Corollary}

\newtheorem{ques}[thm]{Question}

\theoremstyle{definition}
\newtheorem{defn}[thm]{Definition}
\newtheorem{ex}[thm]{Example}
\newtheorem{rmk}[thm]{Remark}

\numberwithin{equation}{section}

\newcommand{\bN}{{\mathbb N}}

\newcommand\Sym{\operatorname{Sym}}

\begin{document}



\title{On soficity for certain fundamental groups of graphs of groups}

\author{David Gao}
\address{\parbox{\linewidth}{Department of Mathematics, University of California, \\
San Diego, 9500 Gilman Drive \# 0112, La Jolla, CA 92093}}
\email{weg002@ucsd.edu}
\urladdr{https://sites.google.com/ucsd.edu/david-gao}
\author{Srivatsav Kunnawalkam Elayavalli}
\address{\parbox{\linewidth}{Department of Mathematics, University of California, \\
San Diego, 9500 Gilman Drive \# 0112, La Jolla, CA 92093}}
\email{skunnawalkamelayaval@ucsd.edu}
\urladdr{https://sites.google.com/view/srivatsavke}
\author{Mahan Mj}
\address{\parbox{\linewidth}{Department of Mathematics, Tata Institute of Fundamental Research, \\
Dr Homi Bhabha Road, Navy Nagar, Colaba, Mumbai, Maharashtra 400005, India
}}
\email{mahan@math.tifr.res.in}
\urladdr{https://mathweb.tifr.res.in/~mahan/}

\thanks{S. Kunnawalkam Elayavalli gratefully acknowledges support from NSF grant DMS 2350049. Mahan Mj is partly supported by a DST JC Bose Fellowship,  by  the Department of Atomic Energy, Government of India, under project no.12-R\&D-TFR-5.01-0500, and by an endowment of the Infosys Foundation.}

\begin{abstract}
In this note we study a family of  graphs of groups over arbitrary base graphs where 
\begin{enumerate}
\item all vertex groups are isomorphic to a fixed countable sofic group $G$,
\item all edge groups $H<G$ are such that the embeddings of $H$ into $G$ are identical everywhere. 
\end{enumerate}
We prove soficity for this family of groups under a flexible technical hypothesis for $H$ called $\sigma$-co-sofic.
This proves soficity for group doubles $*_H G$, where $H<G$ is an arbitrary separable subgroup and $G$ is countable and sofic. This includes arbitrary finite index group doubles of sofic groups among various other examples. 
\end{abstract}

\maketitle

    \section{Introduction}

The notion of soficity for countable groups introduced by Gromov has been very useful to study both algebraic and analytic questions in group theory. For instance, see \cite{surveysofic, KLi, elekszabodirect, Bowenfinvariant, Luck, ThomDiophantine, Hayes5, jung2016ranktheoreml2invariantsfree} where important questions are answered for families of sofic groups. It is therefore  interesting to prove soficity for new families of groups. It is known that soficity is preserved under certain group constructions such as cartesian and free products, amalgamated free products over amenable amalgams \cite{ESZ2}, graph products \cite{CHR}, wreath products \cite{HayesSale}, and some generalized wreath products \cite{gao2024soficitygroupactionssets}. Despite these results, the question of arbitrary amalgamated free products is very difficult to address, a notable example being the Burger-Mozes group \cite{BMgroup}, which arises as an amalgamated free product of two free groups over a certain fixed finite index subgroup. The primary motivation of this note is to address \emph{doubles of groups over subgroups}, i.e.\ the situation of amalgamated free products in which both embeddings of the amalgamating subgroup are identical up to automorphism. 

Let $G$ be a countable group, and $H<G$ a fixed subgroup. Let $\Gamma$ be an arbitrary countable graph. Denote by $D_{\Gamma}(G, H)$ the fundamental group of graph of groups over $\Gamma$ with vertex groups $G$ and edge groups $H$ with all embeddings identical to the original one. To set up the main result, we need the following  definition: 

\begin{defn}
    Let $H < G$ be a pair of groups. We say $H$ is \emph{co-sofic}
    in $G$ if there exist two decreasing sequences $(G_i)_{i \in \bN}$ and $(H_i)_{i \in \bN}$ of subgroups of $G$ s.t. $\cap_i G_i = H$, $H_i < G_i$, $H_i \triangleleft G$, $G/H_i$ is sofic, and $G_i/H_i$ is amenable for all $i$. Now, let $H < G$ be a pair of groups. We say $H$ is $\sigma$\emph{-co-sofic} in $G$ if there exist two increasing sequences $(G_i)_{i \in \bN}$ and $(H_i)_{i \in \bN}$ of subgroups of $G$ s.t. $\cup_i G_i = G$, $\cup_i H_i = H$, $H_i < G_i$, and $H_i$ is co-sofic in $G_i$ for all $i$.
\end{defn} 

We can now state our main result: 

\begin{thm}\label{main theorem}
Let $G$ be a countable sofic groups, and $H<G$ a fixed $\sigma$-co-sofic subgroup. Let $\Gamma$ be an arbitrary countable graph. Then $D_{\Gamma}(G,H)$ is sofic. 
\end{thm}

Now we address some concrete examples arising from our main result. Denote by $\ast_HG$ the fundamental group of graph of groups over the by infinite line graph (Cayley graph of $\mathbb{Z}$)  with vertex groups $G$ and edge groups $H$ with all embeddings being  identical copies of $H<G$.   Say that a subgroup $H<G$ is separable if there exists a sequence of finite index subgroups $H_i<G$ such that $\cap_i H_i=H$. Note that is proved in Example \ref{co-sofic-ex} that separable subgroups are $\sigma$-co-sofic. Therefore we get the following result:

\begin{cor}\label{main theorem A}
    If $G$ is a countable sofic group and $H<G$ is separable, then $\ast_HG$ is sofic. In particular all finite index doubles of sofic groups are sofic.
\end{cor}

 It is a folklore fact, using a normal form argument, that doubles of residually finite groups over separable subgroups are residually finite. Therefore, arbitrary doubles of locally extended residually finite (LERF) groups over finitely generated subgroups are residually finite, hence sofic (this is also a consequence of our result above). (Recall that a group is called LERF if all of its finitely generated subgroups $H$ are separable.) Note that Theorem \ref{main theorem A} can be considered as a sofic analogue of the fact that doubles of residually finite groups over separable subgroups  are residually finite. However, unlike the residually finite case, the proof here requires  more care in handling approximate homomorphisms. 
 
 \begin{rmk}
     A list of examples of LERF groups include: finitely generated nilpotent groups, finitely generated free groups \cite{hall}, surface groups \cite{scott,scott-corr}, limit groups \cite{wiltongafa}, graphs of free groups with cyclic edge groups \cite{brunner, wisepaper}, several Artin groups \cite{ALMEIDA202125}, fundamental groups of compact hyperbolic 3-manifolds \cite{agoldoc}, etc. By passing to normal cores and using a limiting argument, it is  easy to see that arbitrary subgroups of LERF groups are $\sigma$-co-sofic. This is explained in Example \ref{sigma-co-sofic-ex}. Therefore one obtains that all doubles of LERF groups over arbitrary subgroups are sofic. However, in contrast, note that it is not true that all doubles of LERF groups are residually finite. Indeed, one can consider doubles over non finitely generated non separable subgroups, such as $N \triangleleft F$ a normal subgroup of a free group $F$ such that $F/N$ is not residually finite. However, a combination of the fact that an arbitrary double of a LERF group is a limit (in the space of marked groups) of doubles over finitely generated subgroups, the fact that these groups are residually finite, and the fact that soficity is preserved under limits, also implies that these groups are sofic. 
\end{rmk}

\begin{rmk}
    We also remark that our main result can be extended to fundamental groups of graphs of groups over finite graphs with a fixed $G$ as vertex groups, but with different edge groups $H_e$ embedded in the same way in $G$ at the source and end of the edge $e$. In this case, one  defines a suitable analogue of $\sigma$-co-soficity for the family $H_e<G$, in order to prove soficity. We omit substantiating this further  in this note due to technicalities and lack of substantial new examples. 
\end{rmk}
 The proofs of our results are motivated by work of the first author \cite{gaoprepr} on  related topics in the von Neumann algebra setting.  Corollary \ref{main theorem A} follows immediately from  by Theorem \ref{main theorem}, see Example \ref{co-sofic-ex}. The proof of Theorem \ref{main theorem} is presented in the end of Section 2. In conclusion, we leave the following as a natural question arising from this work: 

\begin{ques}
    If $G$ is a countable sofic group, then are arbitrary group doubles $\ast_HG$ are sofic?
\end{ques}


    

    \subsection*{Acknowledgments:} We sincerely thank Marco Linton for helpful comments on references regarding LERF groups and Francesco Fournier-Facio for several helpful comments on a previous draft. 
    \section{Proofs of results}

\begin{defn}
Let $G$ be a countable discrete group, $A$ be a finite set, $F \subseteq G$ be a finite subset, $\epsilon > 0$, and $\varphi: G \rightarrow \Sym(A)$ be a map (not necessarily a homomorphism).

\begin{enumerate}
    \item $\varphi$ is called \emph{unital} if it sends the identity of $G$ to the identity of $\Sym(A)$;
    \item $\varphi$ is called $(F, \epsilon)$\emph{-multiplicative} if $d(\varphi(gh),\varphi(g)\varphi(h)) < \epsilon$ for all $g, h \in F$, where $d$ denotes the normalized Hamming distance;
    \item $\varphi$ is called $(F, \epsilon)$\emph{-free} if $d(\varphi(g),\text{Id}_A) > 1 - \epsilon$ for all non-identity $g \in F$, where $d$ denotes the normalized Hamming distance;
    \item $\varphi$ is called an $(F, \epsilon)$\emph{-approximation} if it is unital, $(F, \epsilon)$-multiplicative, and $(F, \epsilon)$-free;
    \item $G$ is called \emph{sofic} if, for every finite $F \subseteq G$ and $\epsilon > 0$, there exists an $(F, \epsilon)$-approximation.
\end{enumerate}
\end{defn}

Let $\omega$ be a free ultrafilter on $\bN$ and let $H_i$ be a sequence of countable groups. Denote by $\prod_\omega H_i$ the algebraic ultrapower, i.e, $(\times_{i\in \bN} H_i)/N$ where $N=\{(g_i)| \{i: g_i=1\}\in \omega\}$. 
\begin{defn}
Let $H < G$ be a pair of groups, $K$ be a group. We say $H < G$ is \emph{relatively sofic over }$K$ if there exists a sequence of inclusions of groups $(H_i < G_i)_{i \in \bN}$ where all $G_i$ are sofic, all $H_i$ are amenable, and there exists an embedding $\pi: G \rightarrow \prod_\omega (G_i \times K)$ where $\omega$ is a free ultrafilter on $\bN$ and the ultraproduct is algebraic, s.t. $\pi(G) \cap \prod_\omega (H_i \times K) = \pi(H)$.
\end{defn}

Note that $\{e\} < H$ is relatively sofic over a sofic group $K$ is equivalent to $H$ being sofic. This follows from Lemma \ref{soficity-from-ulatrprod}. 

\begin{lemm}\label{commut-sq-embed}
Let $H, K < G$ be inclusions of groups. Then $\ast_{H \cap K} H$ embeds into $\ast_K G$.
\end{lemm}
\begin{proof} We need to define an embedding $\varphi: \ast_{H \cap K} H \rightarrow \ast_K G$. For clarify, we index the free multipliers in $\ast_{H \cap K} H$ and in $\ast_K G$ by an index set $I$, i.e., $H_i$ are copies of $H$ and $G_i$ are copies of $G$ for all $i \in I$. An element of $\ast_{H \cap K} H$ is then either an element of $H \cap K$, or can be written in standard form $h_{1, i_1} \cdots h_{n, i_n}$ where $h_{j, i_j} \in H_{i_j} \setminus (H \cap K) = H_{i_j} \setminus K$ for all $j \in \{1, \cdots, n\}$ and $i_j \neq i_{j + 1}$ for all $j \in \{1, \cdots, n - 1\}$. In the first case, define $\varphi$ to act as the natural embedding of $H \cap K$ into $K$ then into $\ast_K G$. In the second case, define $\varphi(h_{1, i_1} \cdots h_{n, i_n}) = h_{1, i_1} \cdots h_{n, i_n}$. Since $h_{j, i_j} \in H_{i_j} \setminus K \subseteq G_{i_j} \setminus K$, $h_{1, i_1} \cdots h_{n, i_n}$ is already in standard form in $\ast_K G$. It is then easy to show that $\varphi$, thus defined, is a well-defined injective group homomorphism.
\end{proof}

\begin{lemm}\label{prod-embed}
Let $H < G$ be a pair of groups, $K$ be a group. Then $\ast_{H \times K} (G \times K)$ embeds into $(\ast_H G) \times K$.
\end{lemm}
\begin{proof} We need to define an embedding $\varphi: \ast_{H \times K} (G \times K) \rightarrow (\ast_H G) \times K$. For clarify, we index the free multipliers in $\ast_{H \times K} (G \times K)$ and in $\ast_H G$ by an index set $I$, i.e., $(G \times K)_i$ are copies of $G \times K$ and $G_i$ are copies of $G$ for all $i \in I$. An element of $\ast_{H \times K} (G \times K)$ is then either an element of $H \times K$, or can be written in standard form $(g_{1, i_1}, k_{1, i_1}) \cdots (g_{n, i_n}, k_{n, i_n})$ where $(g_{j, i_j}, k_{j, i_j}) \in (G \times K)_{i_j} \setminus (H \times K) = (G_{i_j} \setminus H) \times K$ for all $j \in \{1, \cdots, n\}$ and $i_j \neq i_{j + 1}$ for all $j \in \{1, \cdots, n - 1\}$. In particular, $g_{j, i_j} \in G_{i_j} \setminus H$ for all $j$. In the first case, define $\varphi$ to act as the natural embedding of $H \times K$ into $(\ast_H G) \times K$. In the second case, define $\varphi((g_{1, i_1}, k_{1, i_1}) \cdots (g_{n, i_n}, k_{n, i_n})) = (g_{1, i_1} \cdots g_{n, i_n}, k_{1, i_1} \cdots k_{n, i_n}) \in (\ast_H G) \times K$. Since $g_{j, i_j} \in G_{i_j} \setminus H$, $g_{1, i_1} \cdots g_{n, i_n}$ is in standard form in $\ast_H G$. It is then easy to show that $\varphi$, thus defined, is a well-defined injective group homomorphism.
\end{proof}

\begin{lemm}\label{ultraprod-embed}
Let $(H_i < G_i)_{i \in \bN}$ be a sequence of inclusions of groups. Let $\omega$ be a free ultrafilter on $\bN$ and $\mathbf{G} = \prod_\omega G_i$, $\mathbf{H} = \prod_\omega H_i$ be algebraic ultraproducts. Then $\ast_{\mathbf{H}} \mathbf{G}$ naturally embeds into $\mathbf{K} = \prod_\omega \ast_{H_i} G_i$, where the ultraproduct is again algebraic.
\end{lemm}
\begin{proof} Observe that if $\{i \in \bN: H_i = G_i\} \in \omega$, then $\mathbf{G} = \mathbf{H} = \mathbf{K}$, whence the result is trivial. Thus, we may assume $\{i \in \bN: H_i = G_i\} \notin \omega$. Note that changing $G_i$ and $H_i$ for those $i$ within the set $\{i \in \bN: H_i = G_i\}$ does not change $\mathbf{G}$, $\mathbf{H}$, or $\mathbf{K}$. Hence, by doing so we may assume $H_i \subsetneq G_i$ for all $i$. Now, for any $g \in \mathbf{G} \setminus \mathbf{H}$, we may lift it to a sequence of element $(g_i) \in \prod_{i \in \bN} G_i$. Since $g \notin \mathbf{H}$, $\{i \in \bN: g_i \in H_i\} \notin \omega$, so by replacing all $g_i$ for $i$ within the set $\{i \in \bN: g_i \in H_i\}$, we may assume $g_i \notin H_i$ for all $i$. Call such a sequence $(g_i)$ a proper representing sequence of $g$. We now define a map $\varphi: \ast_{\mathbf{H}} \mathbf{G} \rightarrow \mathbf{K}$. For clarify, we index the free multipliers in $\ast_{\mathbf{H}} \mathbf{G}$ by an index set $I$, i.e., $\mathbf{G}_i$ are copies of $\mathbf{G}$ for all $i \in I$. An element in $\ast_{\mathbf{H}} \mathbf{G}$ is then either an element of $\mathbf{H}$, or can be written in standard form $g_{1, i_1} \cdots g_{n, i_n}$ where $g_{j, i_j} \in \mathbf{G}_{i_j} \setminus \mathbf{H}$ for all $j \in \{1, \cdots, n\}$ and $i_j \neq i_{j + 1}$ for all $j \in \{1, \cdots, n - 1\}$. In the first case, define $\varphi$ to act as the natural embedding of $\mathbf{H} = \prod_\omega H_i$ into $\mathbf{K} = \prod_\omega \ast_{H_i} G_i$. In the second case, for each $g_{j, i_j}$, choose a proper representing sequence $(g_{j, i_j, k})_{k \in \bN}$. Then define $\varphi(g_{1, i_1} \cdots g_{n, i_n})$ to be the element of $\mathbf{K}$ represented by the sequence $(g_{1, i_1, k} \cdots g_{n, i_n, k})_{k \in \bN}$. Observe that $g_{1, i_1, k} \cdots g_{n, i_n, k}$ is an element in standard form in $\ast_{H_i} G_i$, from which it is easy to show that $\varphi$, thus defined, is a well-defined injective group homomorphism.
\end{proof}

\begin{lemm}\label{soficity-from-ulatrprod}
Let $G$ be a countable discrete group, $(H_i)_{i \in \bN}$ be a sequence of sofic groups, $\omega$ be a free ultrafilter on $\bN$. Suppose $G$ embeds into $\prod_\omega H_i$ where the ultraproduct is algebraic, then $G$ is sofic.
\end{lemm}
\begin{proof} For any finite $F \subseteq G$ and $\epsilon > 0$, we need to construct an $(F, \epsilon)$-approximation. Assume WLOG that $F$ contains the identity. Identify $G$ with its image under the embedding into $\prod_\omega H_i$. Then, for each $g \in F \cdot F$, we may lift it to a representing sequence $(g_i)$. Since $F$ is finite, we may choose a set $R \in \omega$ s.t.,

\begin{enumerate}
    \item For each $g, h \in F$, $i \in R$, we have $g_ih_i = (gh)_i$;
    \item For each non-identity $g \in F$, $i \in R$, we have $g_i \neq e$;
    \item For $g = e$, whenever $i \in R$, we have $g_i = e$.
\end{enumerate}

Fix any $i_0 \in R$. Since $H_{i_0}$ is sofic, there exists an $(F_0, \epsilon)$-approximation $\varphi_0: H_{i_0} \rightarrow \Sym(A)$, where $F_0 = \{g_{i_0}: g \in F\}$. We may then define $\varphi: G \rightarrow \Sym(A)$ by $\varphi(g) = \varphi_0(g_i)$ whenever $g \in F \cdot F$ and where $\varphi(g)$ can be defined arbitrarily for $g \notin F \cdot F$. It is easy to verify $\varphi$ is an $(F, \epsilon)$-approximation.
\end{proof}

\begin{thm}\label{thm-rel-soficity}
Suppose $H < G$ is relatively sofic over a sofic group $K$, then $\ast_H G$ is sofic.
\end{thm}
\begin{proof} By definition, there exists a sequence of inclusions of groups $(H_i < G_i)_{i \in \bN}$ where all $G_i$ are sofic, all $H_i$ are amenable, and there exists an embedding $\pi: G \rightarrow \prod_\omega (G_i \times K)$ where $\omega$ is a free ultrafilter on $\bN$, s.t. $\pi(G) \cap \prod_\omega (H_i \times K) = \pi(H)$. By Lemma \ref{commut-sq-embed}, $\ast_H G = \ast_{\pi(H)} \pi(G)$ embeds into $\ast_{\prod_\omega (H_i \times K)} (\prod_\omega (G_i \times K))$, which, by Lemma \ref{ultraprod-embed}, embeds into $\prod_\omega \ast_{H_i \times K} (G_i \times K)$, which, by Lemma \ref{prod-embed}, then embeds into $\prod_\omega ((\ast_{H_i} G_i) \times K)$. By the main result of \cite{ESZ2}, as $H_i$ are amenable and $G_i$ are sofic, $\ast_{H_i} G_i$ are sofic for all $i$. Hence, so are $(\ast_{H_i} G_i) \times K$ as $K$ is sofic. The result then follows by Lemma \ref{soficity-from-ulatrprod}.
\end{proof}

We also note that relative soficity is preserved under finite products. For this, we first need the following elementary fact about ultrafilters:

\begin{lemm}\label{ultra-lemm}
Let $\omega_k$ be a sequence of free ultrafilters on $\bN$ and $\omega$ be a free ultrafilter on $\bN$. Let,
\begin{equation*}
    \varepsilon = \{A \subseteq \bN \times \bN: \{k \in \bN: p_2((\{k\} \times \bN) \cap A) \in \omega_k\} \in \omega\}
\end{equation*}
where $p_2: \bN \times \bN \rightarrow \bN$ is the projection onto the second component. Then $\varepsilon$ is a free ultrafilter on $\bN \times \bN$.
\end{lemm}

\begin{proof}
We first observe that, as $\bN \in \omega_k$ and $\bN \in \omega$, we have $\bN \times \bN \in \varepsilon$. Now, assume $A \in \varepsilon$ and $A \subseteq B$. Then for each $k$, as $\omega_k$ is upward closed, we have $p_2((\{k\} \times \bN) \cap A) \in \omega_k$ implies $p_2((\{k\} \times \bN) \cap B) \in \omega_k$, i.e., $\{k \in \bN: p_2((\{k\} \times \bN) \cap A) \in \omega_k\} \subseteq \{k \in \bN: p_2((\{k\} \times \bN) \cap B) \in \omega_k\}$. As the former belongs to $\omega$ and $\omega$ is upward closed, the latter set belongs to $\omega$ as well, i.e., $B \in \varepsilon$.

Next, let $A, B \subseteq \varepsilon$. Then,
\begin{equation*}
    I = \{k \in \bN: p_2((\{k\} \times \bN) \cap A) \in \omega_k\} \cap \{k \in \bN: p_2((\{k\} \times \bN) \cap B) \in \omega_k\} \in \omega
\end{equation*}
as both sides of the set intersection belong to $\omega$. For any $k \in I$, we have, since $p_2$ is bijective on $\{k\} \times \bN$,
\begin{equation*}
    p_2((\{k\} \times \bN) \cap (A \cap B)) = p_2((\{k\} \times \bN) \cap A) \cap p_2((\{k\} \times \bN) \cap B) \in \omega_k
\end{equation*}

Thus,
\begin{equation*}
    I \subseteq \{k \in \bN: p_2((\{k\} \times \bN) \cap (A \cap B)) \in \omega_k\}
\end{equation*}

Since $I \in \omega$ and $\omega$ is upward closed, the RHS belongs to $\omega$, i.e., $A \cap B \in \varepsilon$.

Next, let $A \subseteq \bN \times \bN$ s.t. $A \notin \varepsilon$. We need to show $A^c \in \varepsilon$. For each $k$, again as $p_2$ is bijective on $\{k\} \times \bN$, we have,
\begin{equation*}
    p_2((\{k\} \times \bN) \cap A^c) = p_2((\{k\} \times \bN) \cap A)^c
\end{equation*}

As $\omega_k$ is an ultrafilter, this implies $p_2((\{k\} \times \bN) \cap A^c) \in \omega_k$ iff $p_2((\{k\} \times \bN) \cap A) \notin \omega_k$, i.e.,
\begin{equation*}
    \{k \in \bN: p_2((\{k\} \times \bN) \cap A^c) \in \omega_k\} = \{k \in \bN: p_2((\{k\} \times \bN) \cap A) \in \omega_k\}^c
\end{equation*}

As $A \notin \varepsilon$, the complement of the latter set is not in $\omega$. As $\omega$ is an ultrafilter, this implies the former set is in $\omega$, i.e., $A^c \in \varepsilon$. This shows $\varepsilon$ is an ultrafilter.

Finally, we show it is free by showing that, if $A \subseteq \bN \times \bN$ is finite, then $A \notin \varepsilon$. Indeed, $A$ being finite implies $p_2((\{k\} \times \bN) \cap A)$ is finite. Since $\omega_k$ is free, $p_2((\{k\} \times \bN) \cap A) \notin \omega_k$ for all $k$. Thus,
\begin{equation*}
    \{k \in \bN: p_2((\{k\} \times \bN) \cap A) \in \omega_k\} = \varnothing \notin \omega
\end{equation*}
so $A \notin \varepsilon$.
\end{proof}

\begin{prop}\label{perm-of-rel-soficity}
Suppose $H_1 < G_1$ and $H_2 < G_2$ are relatively sofic over $K_1, K_2$, respectively, then $H_1 \times H_2 < G_1 \times G_2$ is relatively sofic over $K_1 \times K_2$
\end{prop}

\begin{proof}
By definition, for each $k = 1, 2$, there exists a free ultrafilter $\omega_k$ on $\bN$, a sequence of inclusions of groups $(H_{k,i} < G_{k,i})_{i \in \bN}$ where all $G_{k,i}$ are sofic and all $H_{k,i}$ are amenable, and an embedding $\pi_k: G_k \to \prod_{\omega_k} (G_{k,i} \times K_k)$ s.t. $\pi_k(H_k) = \pi_k(G_k) \cap \prod_{\omega_k} (H_{i,k} \times K_k)$. For each $g \in G_k$, write $\pi_k(g) = (\widetilde{g_i})_{i \to \omega_k}$. Then it is easy to check that,
\begin{equation*}
    \pi: G_1 \times G_2 \to \prod_{i \to \omega_1} \prod_{j \to \omega_2} [(G_{1,i} \times G_{2,j}) \times (K_1 \times K_2)] = \prod_{i \to \omega_1} \prod_{j \to \omega_2} [(G_{1,i} \times K_1) \times (G_{2,j} \times K_2)]
\end{equation*}
defined by,
\begin{equation*}
    \pi(g^1, g^2) = (\widetilde{g^1_i}, \widetilde{g^2_j})_{i \to \omega_1, j \to \omega_2}
\end{equation*}
is a group embedding. Furthermore, $\pi(H_1 \times H_2) = \pi(G_1 \times G_2) \cap \prod_{i \to \omega_1} \prod_{j \to \omega_2} [(H_{1,i} \times H_{2,j}) \times (K_1 \times K_2)]$. We also observer that, as all $G_{k,i}$ are sofic, all $G_{1,i} \times G_{2,j}$ are sofic. And, as all $H_{k,i}$ are amenable, all $H_{1,i} \times H_{2,j}$ are amenable. Finally, observe that $\prod_{i \to \omega_1} \prod_{j \to \omega_2} [(G_{1,i} \times G_{2,j}) \times (K_1 \times K_2)]$ can be naturally identified with the ultraproduct of $(G_{1,i} \times G_{2,j}) \times (K_1 \times K_2)$ with respect to $\varepsilon = \{A \subseteq \bN \times \bN: \{i \in \bN: p_2((\{i\} \times \bN) \cap A) \in \omega_2\} \in \omega_1\}$, where $p_2: \bN \times \bN \rightarrow \bN$ is the projection onto the second component, which is a free ultrafilter on $\bN \times \bN$ per Lemma \ref{ultra-lemm}. This proves the proposition.
\end{proof}

\begin{defn}
Let $G_k$ be a sequence of groups all sharing a subgroup $H$, $K$ be a group. We say the sequence $(H < G_k)$ is \emph{consistently relatively sofic over }$K$ if there exists a sequence of inclusions of groups $(H_i < L_i)_{i \in \bN}$ where all $L_i$ are sofic, all $H_i$ are amenable, and there exists a sequence of embeddings $\pi_k: G_k \rightarrow \prod_\omega (L_i \times K)$ where $\omega$ is a free ultrafilter on $\bN$, s.t.
\begin{enumerate}
    \item $\pi_k(G_k) \cap \prod_\omega (H_i \times K) = \pi_k(H)$ for all $k$;
    \item For all $h \in H$, $\pi_k(h) \in \prod_\omega (L_i \times K)$ is independent of $k$.
\end{enumerate}
\end{defn}

By essentially the same argument for Theorem \ref{thm-rel-soficity}, we have:

\begin{thm}
Suppose $(H < G_k)_{k \in \bN}$ is consistently relatively sofic over a sofic group $K$, then $\ast^{k \in \bN}_H G_k$ is sofic.
\end{thm}

\begin{defn}
Let $H < G$ be a pair of groups. We say $H$ is \emph{co-sofic} in $G$ if there exist two decreasing sequences $(G_i)_{i \in \bN}$ and $(H_i)_{i \in \bN}$ of subgroups of $G$ s.t. $\cap_i G_i = H$, $H_i < G_i$, $H_i \triangleleft G$, $G/H_i$ is sofic, and $G_i/H_i$ is amenable for all $i$.
\end{defn}

\begin{ex}\label{co-sofic-ex}
A finite index subgroup $H < G$ is co-sofic. This is due to the fact that when $[G: H] < \infty$, we have the normal core $H_G = \cap_{g \in G} gHg^{-1}$ is a normal subgroup of $G$ contained in $H$ and of finite index in $G$ itself. Indeed, if $[G: H] = n$ and $G/H = \{\widetilde{g_1}H, \cdots, \widetilde{g_n}H\}$ for $\widetilde{g_1}, \cdots, \widetilde{g_n} \in G$, then $H_G = \cap_{g \in G} gHg^{-1} = \cap_{k=1}^n \widetilde{g_k}H\widetilde{g_k}^{-1}$ is a finite intersection of finite index subgroups, whence finite index itself. Thus, we may simply choose $G_i = H$ and $H_i = H_G$ for all $i$ in the definition of co-soficity. More generally, if $H < G$ is separable, i.e., there exists a decreasing sequence of finite index subgroups $G_i < G$ s.t. $\cap_i G_i = H$, then $H$ is co-sofic in $G$.
\end{ex}

\begin{defn}
Let $H < G$ be a pair of groups. We say $H$ is $\sigma$\emph{-co-sofic} in $G$ if there exist two increasing sequences $(G_i)_{i \in \bN}$ and $(H_i)_{i \in \bN}$ of subgroups of $G$ s.t. $\cup_i G_i = G$, $\cup_i H_i = H$, $H_i < G_i$, and $H_i$ is co-sofic in $G_i$ for all $i$.
\end{defn}

\begin{ex}\label{sigma-co-sofic-ex}
Clearly, if $H < G$ is co-sofic, then it is $\sigma$-co-sofic as well. More generally, if $H < G$ is an increasing union of co-sofic subgroups $H_i < G$, then $H$ is $\sigma$-co-sofic in $G$. Recall that a group $G$ is called \emph{locally extended residually finite (LERF)}    all its finitely generated subgroups are separable. So, by Example \ref{co-sofic-ex}, a LERF $G$ has all its finitely generated subgroups co-sofic in it. Thus, a LERF group $G$ has \emph{all} its subgroups  $\sigma$-co-sofic in it. This is because any subgroup of $G$ is an increasing union of finitely generated subgroups.
\end{ex}

\begin{rmk}\label{rmk-co-sofic-perm}
We observe here that $H_1 < G_1$ and $H_2 < G_2$ being co-sofic implies $H_1 \times H_2 < G_1 \times G_2$ is co-sofic. Indeed, by definition, for each $k = 1, 2$, there exist two decreasing sequences $(G_{k,i})_{i \in \bN}$ and $(H_{k,i})_{i \in \bN}$ of subgroups of $G_k$ s.t. $\cap_i G_{k,i} = H_k$, $H_{k,i} < G_{k,i}$, $H_{k,i} \triangleleft G_k$, $G_k/H_{k,i}$ is sofic, and $G_{k,i}/H_{k,i}$ is amenable for all $i$. But then $(G_{1,i} \times G_{2,i})_{i \in \bN}$ and $(H_{1,i} \times H_{2,i})_{i \in \bN}$ are two decreasing sequences of subgroups of $G_1 \times G_2$, $\cap_i (G_{1,i} \times G_{2,i}) = H_1 \times H_2$, $H_{1,i} \times H_{2,i} < G_{1,i} \times G_{2,i}$, $H_{1,i} \times H_{2,i} \triangleleft G_1 \times G_2$, $(G_1 \times G_2)/(H_{1,i} \times H_{2,i}) = (G_1/H_{1,i}) \times (G_2/H_{2,i})$ is a product of two sofic groups and therefore sofic, and $(G_{1,i} \times G_{2,i})/(H_{1,i} \times H_{2,i}) = (G_{1,i}/H_{1,i}) \times (G_{2,i}/H_{2,i})$ is a product of two amenable groups and therefore amenable. Thus, $H_1 \times H_2 < G_1 \times G_2$ is co-sofic.

As a corollary, we also have $H_1 < G_1$ and $H_2 < G_2$ being $\sigma$-co-sofic implies $H_1 \times H_2 < G_1 \times G_2$ is $\sigma$-co-sofic. The proof is simply writing out the definition. We therefore omit the details.
\end{rmk}

\begin{thm}\label{co-sofic-thm}
Let $H < G$, $H$ co-sofic in $G$. Then $H < G$ is relatively sofic over $G$.
\end{thm}
\begin{proof} By definition, we may choose decreasing sequences $(G_k)_{k \in \bN}$ and $(H_k)_{k \in \bN}$ of subgroups of $G$ s.t. $\cap_k G_k = H$, $H_k < G_k$, $H_k \triangleleft G$, $G/H_k$ is sofic, and $G_k/H_k$ is amenable for all $k$. It now suffices to construct an embedding $\varphi: G \rightarrow \prod_\omega (G/H_k \times G)$ s.t. $\varphi(G) \cap \prod_\omega (G_k/H_k \times G) = \varphi(H)$. It is easy to verify that letting $\varphi(g)$ be the element of $\prod_\omega (G_k/H_k \times G)$ defined by the sequence $((gH_k, g))_{k \in \bN}$ provides the desired map, as for any fixed $g \in G \setminus H$, it must be outside of $G_k$ for sufficiently large $k$.
\end{proof}

\begin{lemm}\label{increasing-union-rel-sofic}
Let $H < G$ be a pair of groups, $(G_k)_{k \in \bN}$ and $(H_k)_{k \in \bN}$ be two increasing sequences of subgroups of $G$ s.t. $\cup_k G_k = G$, $\cup_k H_k = H$. Suppose $H_k < G_k$ are relatively sofic over a fixed group $K$ for all $k$, then $H < G$ is relatively sofic over $K$.
\end{lemm}
\begin{proof} For each $k$, by definition there exists a sequence of inclusions of groups $(H_{k, i} < G_{k, i})_{i \in \bN}$ where all $G_{k, i}$ are sofic, all $H_{k, i}$ are amenable, and there exists an embedding $\pi_k: G_k \rightarrow \prod_{i \rightarrow \omega_k} (G_{k, i} \times K)$ where $\omega_k$ is a free ultrafilter on $\bN$, s.t. $\pi_k(G_k) \cap \prod_{i \rightarrow \omega_k} (H_{k, i} \times K) = \pi_k(H_k)$. Fix any free ultrafilter $\omega$ on $\bN$. We define a map $\varphi: G \rightarrow \prod_{k \rightarrow \omega} \prod_{i \rightarrow \omega_k} (G_{k, i} \times K)$ by letting $\varphi(g)$ be the element represented by the sequence $(\pi_k(g'_k))_{k \in \bN}$ where $g'_k = e$ if $g \notin G_k$ and $g'_k = g$ otherwise. Noting that, as $G_k$ increases to $G$, for any fixed $g \in G$ it is contained in $G_k$ for all sufficiently large $k$. From this it is easy to verify that $\varphi$ is a well-defined injective group homomorphism. It is also easy to verify that $\varphi(H) = \varphi(G) \cap \prod_{k \rightarrow \omega} \prod_{i \rightarrow \omega_k} (H_{k, i} \times K)$. Observe that $\prod_{k \rightarrow \omega} \prod_{i \rightarrow \omega_k} (G_{k, i} \times K)$ can be naturally identified with the ultraproduct of $G_{k, i} \times K$ with respect to $\varepsilon = \{A \subseteq \bN \times \bN: \{k \in \bN: p_2((\{k\} \times \bN) \cap A) \in \omega_k\} \in \omega\}$, where $p_2: \bN \times \bN \rightarrow \bN$ is the projection onto the second component, which is a free ultrafilter on $\bN \times \bN$ per Lemma \ref{ultra-lemm}. This proves the lemma.
\end{proof}

\begin{cor}\label{sigma-co-sofic-to-rel-sofic}
Let $H < G$, $H$ $\sigma$-co-sofic in $G$. Then $H < G$ is relatively sofic over $G$.
\end{cor}
\begin{proof} This follows immediately from Theorem \ref{co-sofic-thm} and Lemma \ref{increasing-union-rel-sofic}.
\end{proof}

\begin{cor}\label{main cor}
Let $H < G$, $H$ $\sigma$-co-sofic in $G$, $G$ being sofic. Then $\ast_H G$ is sofic.
\end{cor}
\begin{proof} This follows immediately from Theorem \ref{thm-rel-soficity} and Corollary \ref{sigma-co-sofic-to-rel-sofic}.
\end{proof}

Per Example \ref{co-sofic-ex}, this implies the group doubles of sofic groups over separable subgroups are sofic. Per Example \ref{sigma-co-sofic-ex}, and as LERF groups are residually finite and therefore sofic, this implies arbitrary group doubles of LERF groups are sofic. 

\begin{proof}[Proof of Theorem \ref{main theorem}]
    Using either Remark \ref{rmk-co-sofic-perm} or Proposition \ref{perm-of-rel-soficity}, and observing that $\{0\} < \mathbb{Z}$ is co-sofic, we see that $H < G \times \mathbb{Z}$ is relatively sofic over $G \times \mathbb{Z}$. As $G$ is sofic, so is $G \times \mathbb{Z}$. Thus, $\ast_H (G \times \mathbb{Z})$ is sofic. By construction of $D_\Gamma(G, H)$, it is an amalgamated free product of copies of $G$ and copies of $H \times \mathbb{Z}$, over $H$, whence it embeds into $\ast_H (G \times \mathbb{Z})$. Hence, $D_\Gamma(G, H)$ is sofic.
\end{proof}
\bibliographystyle{plain}
\bibliography{soficconj}

\end{document}